\documentclass{amsart}

 \newtheorem{thm}{Theorem}[subsection]
 \newtheorem{cor}[thm]{Corollary}
 \newtheorem{lem}[thm]{Lemma}
 \newtheorem{prop}[thm]{Proposition}
 \theoremstyle{definition}
 \newtheorem{defn}[thm]{Definition}
 \theoremstyle{remark}
 
 \numberwithin{equation}{subsection}

 \newcommand{\norm}[1]{\left\Vert#1\right\Vert}
 
 \newcommand{\C}{\mathbb{C}}
 \DeclareMathOperator{\adj}{adj}
\begin{document}

\title[On Jordan models of $C_{0}$-operators in a multiply connected domain]
 {Linear Algebraic Properties for Jordan Models of $C_{0}$-operators
 relative to multiply connected domains}

\author{ Yun-Su Kim }

\address{Department of Mathematics, Indiana University, Bloomington,
Indiana, U.S.A. }

\email{kimys@indiana.edu}


\dedicatory{}

\commby{Daniel J. Rudolph}


\begin{abstract}We study $C_{0}$-operators relative to a multiply connected
domain using a substitute of the characteristic function. This
method allows us to prove certain relations between the Jordan
model of an operator and that of its restriction to an invariant
subspace.
\end{abstract}

\maketitle

\section*{Introduction}
Hasumi \cite{12}, Sarason \cite{17}, and Voivhick \cite{19}
started operator theory related to function theory on multiply
connected domains by providing an analogue (in the scalar case) of
Beurling's theorem on invariant subspaces of the Hardy spaces of
the open unit disk. Their work was continued in the work of
Abrahamse$-$Douglas [1, 2], and of Ball [4, 5]. In particular,
J.A. Ball \cite{4} introduced the class of $C_{0}$-operators
relative to a bounded finitely connected region $\Omega$ in the
complex plane, whose boundary $\partial\Omega$ consists of a
finite number of disjoint, analytic, simple closed curves. J.
Agler [3] showed that the existence of normal boundary dilations
$-$ an analogue of Sz.-Nagy dilation theorem $-$ still holds for
annuli but it may fail for domains of connectivity greater than
two (Dritschel$-$McCullough \cite{11}). However it holds up to
similarity (Douglas$-$Paulsen \cite{10}); this allowed Zucchi
\cite{20} to provide a classification of $C_{0}$-operators
relative to $\Omega$. Since no analogue of the characteristic
function of a contraction is available in that context, that study
does not yield some of the results available for the unit disk. In
this paper we use a substitute for the characteristic function,
suggested by an analogue of Beurling's theorem provided by M.A.
Abrahamse and R.G. Douglas [2]. This allows us to prove a
relationship between the Jordan models of a $C_{0}$-operator
relative to $\Omega$, of its restriction to an invariant subspace,
and of its compression to the orthocomplement of that subspace. In
the case of the open unit disk, this result was proved by H.
Bercovici and D. Voiculescu [7].

This paper is organized as follows. Section 1 contains
preliminaries about bundle shifts and operators of class $C_{0}$.
Here we define the notion of an operator-valued quasi-inner
function and prove a useful reformulation of the description of
invariant subspaces given in [2].

In Section 2, we review concepts relating quasi-equivalence and
quasi-similarity, which were first introduced in [13, 14] and we
prove the main result.

The author would like to express her gratitude to her thesis
advisor, Professor Hari Bercovici.\vskip5cm
\section{Preliminaries and Notation}

\goodbreak
In this paper, $\C$, $\overline{M}$, and $L(K_{1},$ $K_{2})$
denote the set of complex numbers, the (norm) closure of a set
$M$, and the set of bounded linear operators from $K_{1}$ to
$K_{2}$ where $K_{1}$ and $K_{2}$ are Hilbert spaces,
respectively.
  \subsection{Hardy spaces}
We refer to [16] for basic facts about Hardy space, and recall
here the basic definitions.
\begin{defn}
The space ${\: H^{2}(\Omega)}$ is defined to be the space of
analytic functions $f$ on $\Omega$ such that the subharmonic
function $|f|^{2}$ has a harmonic majorant on $\Omega$. For a
fixed $z_{0}$ $\in\Omega$, there is a norm on
$H^{2}$$($$\Omega$$)$ defined by

      $\|f\|$=inf$\{ u(z_{0})^{1/2}$: $u$ is a harmonic majorant of
      $|f|^{2}\}$.
\end{defn}

 Let $m$ be harmonic measure for the point $z_{0}$, let
 $L^{2}(\partial{\Omega})$ be the $L^{2}$-space of complex valued
 functions on the boundary of $\Omega$ defined with respect to $m$,
 and let $H^{2}(\partial{\Omega})$ be the set of
 functions $f$ in $L^{2}(\partial{\Omega})$ such that
 $\int_{\partial{\Omega}} f(z)g(z) dz$ = 0 for every $g$ that is
 analytic in a neighborhood of the closure of $\Omega$.
 If $f$ is in $H^{2}(\Omega)$, then there is a function $f^{\ast}$ in
 $H^{2}(\partial{\Omega})$ such that $f({z})$ approaches
 $f^{\ast}(\lambda_{0})$ as $z$ approaches $\lambda_{0}$
 nontangentially, for almost every $\lambda_{0}$ relative to $m$. The map
 $f\rightarrow{f^{\ast}}$ is an isometry from $H^{2}(\Omega)$ onto
 $H^{2}(\partial{\Omega})$. In this way, $H^{2}(\Omega)$ can be
 viewed as a closed subspace of $L^{2}(\partial{\Omega})$.

 A function $f$ defined on $\Omega$ is in $H^{\infty}(\Omega)$ if
 it is holomorphic and bounded. $H^{\infty}(\Omega)$ is a closed
 subspace of $L^{\infty}({\Omega})$ and it is a Banach
 algebra if endowed with the supremum norm. Finally, the mapping
 $f\rightarrow{f^{\ast}}$ is an isometry of $H^{\infty}(\Omega)$
 onto a week$^{*}$-closed subalgebra of
 $L^{\infty}(\partial{\Omega})$.

\begin{defn} If $K$ is a Hilbert space, then $H^{2}(\Omega$,$K)$ is
defined to be the space of analytic functions
$f:\Omega\rightarrow{K}$ such that the subharmonic function
${\norm{f}}^{2}$ is majorized by a harmonic function $\nu$. Fix a
point $z_{0}$ in $\Omega$ and define a norm on $H^{2}(\Omega$,$K)$
by

$\norm{f}$=inf $\{{\nu(z_{0})}^{1/2}$ : $\nu$ is a harmonic
majorant of ${\norm{f}}^{2}\}$.
\end{defn}

As before, $H^{2}(\Omega$,$K)$ can be identified with a closed
subspace of the space $L^{2}(\partial\Omega$,$K)$ of square
integrable $K$-valued functions on $\partial\Omega$. Define
$S_{K}:H^{2}(\Omega$,$K)$$\rightarrow$$H^{2}(\Omega$,$K)$ by
$(S_{K}f)(z)$=$zf(z)$.

\subsection{Vector Bundles}

We present in this section and in section $1.3$ the standard
definitions of analytic vector and flat unitary vector bundles. We
refer to [2] for this material.

Let $K$ be a Hilbert space. An \emph{analytic vector bundle} over
$\Omega$ with fiber $K$ is a pair $(E,p)$, where
$p:E\rightarrow\Omega$ is a continuous surjective map such that:

$(1)$ Each $z\in\Omega$ has a neighborhood $U_{z}$ for which there
is a homeomorphism $\varphi_{z}:
{U_{z}}\times$${K}\rightarrow{p^{-1}(U_{z})}$ satisfying
$\varphi_{z}(\omega,k)\in{p^{-1}(\omega)}$ for $\omega\in{U_{z}}$
and $k\in{K}$.

$(2)$ If $z_{1}$, $z_{2}\in\Omega$, there is an analytic map
$\psi_{z_{1},z_{2}}:U_{z_{1}}\cap{U_{z_{2}}}\rightarrow{GL(K)}$
satisfying
$\varphi_{z_{1}}(\omega,k)=\varphi_{z_{2}}(\omega,\psi_{z_{1},z_{2}}(\omega)k)$,
where $GL(K)$ is the set of all invertible linear operators on
$K$. \vskip0.2cm If we can choose $U_{z}=\Omega$ for some
$z\in\Omega$, we say that $(E,p)$ is a \emph{trivial bundle}. If
each $\psi_{z_{1},z_{2}}$ is a constant unitary operator for every
$z_{1}$, $z_{2}\in\Omega$, then $(E,p)$ is called a \emph{flat
unitary vector bundle}. \vskip0.1cm \textbf{Theorem A} $[8]$.
\emph{Every analytic vector bundle over $\Omega$ is analytically
trivial.}
\subsection{Bundle Shift}

Let $E$ be a vector bundle over $\Omega$. A \emph{cross section}
of a vector bundle $E$ over $\Omega$ is a continuous function $f$
from $\Omega$ into $E$ such that $p(f(z))=z$ for all $z$ in
$\Omega$. For each $\omega$ in $U_{z}$, define a map
$\varphi_{z}^{\omega}:K\rightarrow{p^{-1}(\omega)}$ by
$\varphi_{z}^{\omega}(k)=\varphi_{z}(\omega,k)$.

If $E$ is a flat unitary vector bundle over $\Omega$ with fiber
$K$ and if $f$ is a cross section of $E$, then for $\omega$ in
$U_{z_{1}}\cap{U_{z_{2}}}$($z_{1}$, $z_{2}\in\Omega$), the
operator
$(\varphi_{z_{1}}^{\omega})^{-1}$$\varphi_{z_{2}}^{\omega}$ is
unitary so that
$\norm{(\varphi_{z_{2}}^{\omega})^{-1}(f(z))}$$=$$\norm{(\varphi_{z_{1}}^{\omega})^{-1}(f(z))}$.
This means that there is a function $h_{f}:\Omega\rightarrow{R}$
defined by
$h_{f}^{E}$$(z)$=$\norm{(\varphi_{z_{2}}^{\omega})^{-1}(f(z))}$,
where $\omega$ is in $U_{z_{2}}$.
\begin{defn} We define ${\: {H}^{2}(\Omega,E)}$ to be the space of
analytic cross sections $f$ of $E$ such that $(h_{f}^{E})^{2}$ is
majorized by a harmonic function.\end{defn}
 We can define
the\emph{ bundle shift} $T_{E}$ on ${H}^{2}$$(\Omega$,$E)$ by
$(T_{E}f)(z)$=$zf(z)$ for $z\in\Omega$. The operator $T_{E}$
admits a functional calculus defined on the algebra $R(\Omega)$ of
rational functions with poles off $\overline{\Omega}$. More
precisely, if $u\in{R(\Omega)}$, $(u(T_{E})f)(z)=u(z)f(z)$ for
$z\in\Omega$ and $f\in{{H}^{2}(\Omega,E)}$. \vskip 0.2 cm

\subsection{Quasi-Inner Function.}
If $E$ and $F$ are flat unitary bundles over $\Omega$ that extend
to an open set  $\Omega^{\prime}$ containing the closure of
$\Omega$, and $\Theta$ is a bounded holomorphic bundle map from
$E$ to $F$, then $\Theta$ can be shown to have nontangential
limits a.e. relative to $m$ on $\partial\Omega$. The limit at a
point $z$ of $\partial\Omega$ can be regarded as an operator from
the fiber of $E$ at $z$ to the fiber of $F$ at $z$.
\begin{defn} (a) A bounded holomorphic bundle map $\Theta$ is \emph{inner }if
the nontangential limits are isometric operators a.e. relative to
$m$.

(b) Let $K$ and $K^{\prime}$ be Hilbert spaces and let
${H}^{\infty}$$(\Omega$,$L(K$,$K^{\prime}))$ be the Banach space
of all analytic functions
$\Phi:\Omega\rightarrow{L(K,K^{\prime})}$ with the supremum norm.
For $\varphi$$\in$${H}^{\infty}$$(\Omega$,$L(K$,$K^{\prime}))$, we
will say that $\varphi$ is \emph{quasi-inner} if there exists a
constant $c>0$ such that for every $k\in{K}$ and almost every
$z\in\partial\Omega$ we have $\norm{\varphi(z)k}\geq{c\norm{k}}$.
\end{defn}\textbf{Theorem B} [2].\emph{ Let $T_{E}$ be
a bundle shift on ${H}^{2}$$(\Omega$,$E)$. Then a closed subspace
$M$ of ${H}^{2}$$(\Omega$,$E)$ is invariant under the algebra
$\{u(T_{E}):u\in{R(\Omega)}\}$ if and only if
$M$=$\Theta$${H}^{2}$$(\Omega$,$F)$, where $F$ is a flat unitary
bundle over $\Omega$ and $\Theta$ is an inner bundle map from $F$
to $E$.} \vskip0.2cm It will be convenient to reformulate Theorem
B in terms of quasi-inner functions without use of vector bundles.
We will say that a space $M$ is $R(\Omega)$-\emph{invariant} for
an operator $T$ if it is invariant under $u(T)$ for every
$u\in{R(\Omega)}$. For a Hilbert space $K$, define an operator
$S_{K}$ on ${H}^{2}$$(\Omega$,$K)$ by $(S_{K}f)(z)$=$zf(z)$ for
$z\in\Omega$. \vskip0.2cm The proper setting here is maps of flat
unitary vector bundles, i.e., multiplicative multivalued
operator-valued functions. We will convert these to usual single
valued analytic functions by composing them with some bundle
isomorphisms. This has been done quite often in the scalar case,
see, e.g., Royden \cite{15}.

\begin{thm} Let $K$ be a Hilbert space. Then a closed subspace
$M$ of ${H}^{2}$$(\Omega$, $K)$ is $R(\Omega)$-invariant for
$S_{K}$ if and only if there is a Hilbert space $K^{\prime}$ and a
quasi-inner function
$\varphi:\Omega\rightarrow$$L(K^{\prime}$,$K)$ such that
$M=\varphi$${H}^{2}$$(\Omega$,$K^{\prime})$.\end{thm}
\begin{proof} It is clear that a subspace of the form $\varphi$${H}^{2}$$(\Omega$,$K^{\prime})$
with $\varphi:\Omega\rightarrow$$L(K^{\prime}$,$K)$ quasi-inner,
is $R(\Omega)$-invariant. Conversely, consider a closed subspace
$M\subset{H^{2}(\Omega,K)}$ which is $R(\Omega)$-invariant. Let
$M^{\prime}=\{$$G\in$${H}^{2}(\Omega,\Omega\times{K})$: $G(z)=
(z,g(z))$ for some $g\in$$M$$\}$. Then $M^{\prime}$ is a closed
subspace of ${H}^{2}(\Omega,\Omega\times{K})$ which is
$R(\Omega)$-invariant for $T_{\Omega\times{K}}$ and so, by Theorem
B, there is a flat unitary bundle $F$ over $\Omega$ with fiber
$K^{\prime}$, and an inner bundle map
$\Theta:{F}\rightarrow\Omega\times{K}$, such that
$M^{\prime}=\Theta$${H}^{2}$$(\Omega,$$F)$. We know that there is
a flat unitary vector bundle $F^{\prime}$ over an open set
$\Omega^{\prime}$ containing the closure of $\Omega$, with fiber
$K^{\prime}$, such that $F$ is unitary equivalent to the bundle
$F^{\prime}$$\mid$${\Omega}$ \cite{2}. By Theorem A, there is an
analytic isomorphism $\Lambda:
\Omega^{\prime}\times{K^{\prime}}\rightarrow{F^{\prime}}$.

Define an invertible operator $W:
H^{2}(\Omega$,$K^{\prime})\rightarrow{H^{2}(\Omega,
F^{\prime}}$$\mid$${\Omega})$ by $(Wf)(z)=
\Lambda(z,f(z))=\Lambda_{z}(f(z))$ for
$f$$\in$$H^{2}(\Omega$,$K^{\prime})$. Then
$M^{\prime}=\Theta{U}{W}H^{2}(\Omega,K^{\prime})$ where $U$:
$H^{2}(\Omega$,$F^{\prime}$$\mid$$\Omega)$$\rightarrow$
$H^{2}(\Omega$,$F)$ is a unitary operator. For each
$z\in{\Omega}$, we can define a bounded operator $W_{z}$:
$K^{\prime}\rightarrow$$F_{z}$ by $W_{z}a$=$(U(Wh_{a}))(z)$ for
$a\in$$K^{\prime}$ where $h_{a}\in$$H^{2}(\Omega,K^{\prime})$
defined by $h_{a}(z)=a$.

Let $\varphi(z)$=$\Theta_{z}W_{z}$ for $z\in\Omega$ where
$\Theta_{z}=\Theta|{F_{z}}$. Then
$\varphi\in$$H^{\infty}(\Omega,L(K^{\prime},K))$ and
$M=\varphi$${H}^{2}$$(\Omega$,$K^{\prime})$. To conclude our
proof, we must verify that $\varphi$ is quasi-inner.

From the fact that $\Lambda$ is an analytic isomorphism, we see
that the function $z\rightarrow(\Lambda_{z})^{-1}$ is holomorphic
on $\Omega^{\prime}$, and so there is $m>{0}$ such that
$\norm{(\Lambda_{z})^{-1}}$$\leq{m}$ for any $z$$\in$$\Omega$.
Therefore $\norm{W_{z}^{-1}}\leq{m}$ for any $z$$\in$$\Omega$ as
well, so that $\norm{a}/m$$\leq\norm{\varphi(z)a}$ a.e. on
$\partial\Omega$ for $a\in{K^{\prime}}$ as desired.
\end{proof}
\begin{lem}
Let $K_{1}$ and $K_{2}$ be separable Hilbert spaces. If
$T:{H}^{2}(\Omega,K_{1})\rightarrow{H}^{2}(\Omega,K_{2})$ is a
bounded linear operator such that $TS_{K_{1}}=S_{K_{2}}T$, then
there is a function $\psi\in$$H^{\infty}(\Omega,L(K_{1},K_{2}))$
such that $T=M_{\psi}$, where $M_{\psi}(g)(z)=\psi(z)g(z)$ for
$g\in$${H}^{2}(\Omega,K_{1})$, we have
$\norm{T}=\norm{\psi}_{\infty}$.
\end{lem}
\begin{proof}
Define $Y\in$ $(S_{{K_{1}}\bigoplus{K_{2}}})^{\prime}$ by
$Y$$=$$\begin{pmatrix} 0&0\\T&0
\end{pmatrix}$. Then by the proposition 1.9 in [2],
$Y=M_{\omega}$ where
$\omega\in$$H^{\infty}(\Omega,L(K_{1}\bigoplus{K_{2}}))$. Let
$\omega=\begin{pmatrix}
\omega_{11}&\omega_{12}\\\omega_{21}&\omega_{22}
\end{pmatrix}$. Take $\psi=\omega_{21}$, then $T=M_{\psi}$ and one
can check easily that $\norm{T}=\norm{\psi}_{\infty}$.
\end{proof}
\begin{cor} Let $\varphi_{1}
:\Omega\rightarrow{L(K_{1},K)}$ and
$\varphi_{2}:\Omega\rightarrow{L(K_{2},K)}$ be quasi-inner
functions. Then the subspaces $\varphi_{1}H^{2}($$\Omega$,$K_{1})$
and $\varphi_{2}H^{2}($$\Omega$,$K_{2})$ of $H^{2}($$\Omega$,$K)$
are equal if and only if there exist functions
 $\varphi$$\in$$H^{\infty}($$\Omega$,$L($$K_{1}$,
$K_{2})$$)$ and
$\psi\in$$H^{\infty}($$\Omega$,$L($$K_{2}$,$K_{1})$$)$ such that
$\varphi\psi=I_{K_{2}}$, $\psi\varphi=I_{K_{1}}$ and
$\varphi_{1}(z)$=$\varphi_{2}(z)\varphi(z)$ for any
$z$$\in$$\Omega$. In particular, $K_{1}$ and $K_{2}$ have the same
dimension.
\end{cor}
\begin{proof}
The condition $\varphi_{1}(z)$=$\varphi_{2}(z)\varphi(z)$ with
$\varphi$ invertible clearly implies
$\varphi_{1}H^{2}($$\Omega$,$K_{1})$
$=$$\varphi_{2}H^{2}($$\Omega$,$K_{2})$. Conversely, assume that
$\varphi_{1}H^{2}($$\Omega$,$K_{1})$$=$$\varphi_{2}H^{2}($$\Omega$,$K_{2})$.
Define an operator
$T:{H^{2}(\Omega,K_{1})}\rightarrow{H^{2}(\Omega,K_{2})}$ as
follows. For $f\in{H^{2}(\Omega,K_{1})}$, $Tf=g$ such that
$\varphi_{1}f=\varphi_{2}g$. Since $\varphi_{i}$ $(i=1,2)$ is a
quasi-inner function, $T$ is well-defined and invertible. Since
$S_{K_{2}}T=TS_{K_{1}}$, by the previous Lemma $T=M_{\varphi}$ for
a function
$\varphi$$\in$$H^{\infty}($$\Omega$,$L($$K_{1}$,$K_{2})$$)$. Note
that the invertibility of $T$ is equivalent to the invertibility
of $\varphi$. It follows that $\varphi_{1}f=\varphi_{2}\varphi{f}$
for any $f\in{H^{2}(\Omega,K_{1})}$ and so
$\varphi_{1}=\varphi_{2}\varphi$. Since $\varphi(z)$ is invertible
for any $z\in{\Omega}$, $K_{1}$ and $K_{2}$ have the same
dimension.
\end{proof}
\subsection{The Class $C_{0}$}
The theory of Jordan models for contractions of class $C_{0}$ was
developed by Sz.-Nagy$-$Foias, Moore$-$Nordgren, and
Bercovici$-$Voiculescu.

We will present in this section the definition of
$C_{0}$-operators relative to $\Omega$. Reference for this
material is Zucchi [20].

Let ${H}$ be a Hilbert space and ${K_{1}}$ be a compact subset of
the complex plane. If $T$$\in$$L(H)$ and
$\sigma(T)$$\subseteq$$K_{1}$, for $r=p/q$ a rational function
with poles off $K_{1}$, we can define an operator $r(T)$ by
$q(T)^{-1}p(T)$.
\begin{defn}
If $T$$\in$$L(H)$ and $\sigma(T)$$\subseteq$$K_{1}$, we say that
$K_{1}$ is a \emph{spectral set} for the operator $T$ if
$\norm{r(T)}$$\leq$$\max\{$$|$$r(z)$$|$$:$ $z$$\in$$K_{1}\}$,
whenever $r$ is a rational function with poles off $K_{1}$.
\end{defn}
If $T$ $\in$$L(H)$ is an operator with $\overline{\Omega}$ as a
spectral set and with no normal summand with spectrum in
$\partial\Omega$, i.e., $T$ has no reducing subspace
$M$$\subseteq$$H$ such that $T|M$ is normal and
$\sigma(T|M)$$\subseteq$$\partial\Omega$, then we say that $T$
satisfies \emph{hypothesis (h)}.
\begin{thm} ([20], Theorem 3.1.4)
Let $T$$\in$$L(H)$ be an operator satisfying ${\:hypothesis}$
$(h)$. Then there is a unique norm continuous representation
$\Psi_{T}$ of $H^{\infty}(\Omega)$ into $L(H)$ such that :

$(i) \Psi_{T}(1)$=$I_{H}$, where $I_{H}$$\in$$L(H)$ is the
identity operator;

$(ii) \Psi_{T}(g)$=$T$, where $g(z)$=$z$ for all $z$$\in$$\Omega$;

$(iii) \Psi_{T}$ is continuous when $H^{\infty}(\Omega)$ and
$L(H)$ are given the $weak^{*}$-topology.

Moreover $\Psi_{T}$ is contractive, i.e.,
$\norm{\Psi_{T}(f)}$$\leq$$\norm{f}$ for all
$f$$\in$$H^{\infty}(\Omega)$.
\end{thm}

From now on we will indicate $\Psi_{T}(f)$ by $f(T)$ for all
$f$$\in$$H^{\infty}(\Omega)$.

\begin{defn}
An operator $T$ satisfying hypothesis (h) is said to be of
\emph{class} $C_{0}$ \emph{relative to} $\Omega$ if there exists
$u$ $\in$ $H^{\infty}(\Omega)$$\setminus$$\{{0}\}$ such that
$u(T)$=$0$.
\end{defn}
By Theorem 1 in \cite{15}, if $T$ is of class $C_{0}$ relative to
$\Omega$, then there is a quasi-inner function
$m_{T}\in{H^{\infty}(\Omega)}$ such that
$\ker(\Psi_{T})=m_{T}{H^{\infty}(\Omega)}$ and $m_{T}$ is said to
be a \emph{minimal function} of $T$.

\subsection{Jordan Model}

\begin{defn} Let $H$ and $H^{\prime}$ be Hilbert spaces.
An operator $T$$\in$$L(H)$ is called a \emph{quasiaffine
transform} of an operator $T^{\prime}$$\in$$L(H^{\prime})$ if
there exists an injective operator $X$$\in$$L(H,H^{\prime})$ with
dense range in $H^{\prime}$ such that $T^{\prime}X$=$XT$. We write
$T\prec$ $T^{\prime}$ if $T$ is a quasiaffine transform of
$T^{\prime}$. The operators $T$ and $T^{\prime}$ are
\emph{quasisimilar} $(T\sim{T^{\prime}})$ if $T\prec{T^{\prime}}$
and $T^{\prime}\prec{T}$.
\end{defn}

Let $\theta$ and $\theta^{\prime}$ be two functions in
$H^{\infty}(\Omega)$. We say that $\theta$ \emph{divides}
$\theta^{\prime}$ $($or $\theta$$\mid$$\theta^{\prime}$$)$ if
$\theta^{\prime}$ can be written as
$\theta^{\prime}$=$\theta$$\cdot$$\phi$ for some
$\phi$$\in$$H^{\infty}(\Omega)$.  We will use the notation
$\theta\equiv{\theta}^{\prime}$ if
$\theta$$\mid$${\theta}^{\prime}$ and
${\theta}^{\prime}$$\mid$${\theta}$.
\begin{defn}
(i) Given a quasi-inner function
$\theta$$\in$$H^{\infty}(\Omega)$, the \emph{Jordan block
}$S(\theta)$ is the operator acting on the space
$H(\theta)$=$H^{2}(\Omega)\ominus\theta$$H^{2}(\Omega)$ as
follows:\begin{equation}
S(\theta)=P_{H(\theta)}S|H(\theta)\end{equation} where
$S$$\in$$L(H^{2}(\Omega))$ is defined by $(Sf)(z)$=$zf(z)$.

(ii) Let
$\Theta$=$\{\theta_{i}$$\in$$H^{\infty}(\Omega)$$:i=1,2,3,\cdot\cdot\cdot\}$
be a family of quasi-inner functions. Then $\Theta$ is called a
\emph{model function} if $\theta_{i}\mid\theta_{j}$ whenever
$j\leq{i}$. The \emph{Jordan operator} $S(\Theta)$ determined by
the model function $\Theta$ is the $C_{0}$-operator defined as
$S(\Theta)$=$\bigoplus_{i<\gamma^{\prime}}$$S(\theta_{i})$,
$\gamma^{\prime}$=min$\{k$$:$ $\theta_{k}\equiv1\}$.
\end{defn}
 We will call $S(\Theta)$ the \emph{Jordan model }of the operator
$T$ if $S(\Theta)\sim{T}$. From [20], we can get following
results: \vskip0.1cm$\textbf{Theorem}$ $\textbf{C}$. \emph{For
every operator $T$ of class $C_{0}$ relative to $\Omega$ acting on
a separable space $H$, there is a unique Jordan model for $T$.}
\begin{prop}
Let $T$ be of class $C_{0}$ relative to $\Omega$ acting on a
separable space $H$ and $H^{\prime}$ be $R(\Omega)$-invariant for
$T$. If $T\sim{\bigoplus_{\alpha<\gamma}S(\theta_{\alpha})}$ and
$T|H^{\prime}\sim{\bigoplus_{\alpha<\gamma^{\prime}}S(\theta^{\prime}_{\alpha})}$,
then $\theta^{\prime}_{\alpha}|\theta_{\alpha}$ for every
$\alpha\leq$ min$\{\gamma,\gamma^{\prime}\}$.
\end{prop}
\goodbreak
\subsection{Scalar Multiples}
Let $K$ and $K^{\prime}$ be Hilbert spaces and
$\varphi$$\in$${H}^{\infty}$$(\Omega$,$L(K$,$K^{\prime}$$))$ be a
quasi-inner function. We set
$H(\varphi)=H^{2}($$\Omega,$$K^{\prime}$$)\ominus\varphi$$H^{2}($$\Omega,$$K$$)$
and denote by $S(\varphi)$ the compression of $S_{K^{\prime}}$ to
$H(\varphi)$, i.e., $S(\varphi)=P_{H(\varphi)}$
$S_{K^{\prime}}$$|$${H(\varphi)}$, where $P_{H(\varphi)}$ denotes
the orthogonal projection onto $H(\varphi)$.
\begin{defn}
The function $\varphi$$\in$${H}^{\infty}(\Omega,L(K,K^{\prime}))$
is said to have a \emph{scalar multiple}
$u\in{{H}^{\infty}}$$(\Omega)$, $u$$\neq$$0$, if there exists
$\psi$$\in$${H}^{\infty}$$(\Omega$,$L(K^{\prime}$,$K$$)$$)$
satisfying the relation $\varphi(z)$$\psi(z)=
 u(z)I_{K^{\prime}}$ for $z\in\Omega$
\end{defn}
\begin{thm}
Suppose that
$\varphi$$\in$${H}^{\infty}$$(\Omega$,$L(K$,$K^{\prime}$$)$$)$ is
a quasi-inner function and $u$$\in$${H}^{\infty}$$(\Omega$$)$.
Then  the following assertions are equivalent : \vskip 0.1cm $(a)$
$u$ is a scalar multiple of $\varphi$.

$(b)$ $u(S(\varphi))=0$.

$(c)$
$uH^{2}($$\Omega,$$K^{\prime}$$)\subset\varphi$$H^{2}($$\Omega,$$K$$)$.
\end{thm}
\begin{proof}
Assume $(a)$, and let
$\psi$$\in$${H}^{\infty}$$(\Omega$,$L(K^{\prime}$,$K$$)$$)$
satisfy the relation $\varphi(z)$$\psi(z)$$=$$u(z)\cdot$
$I_{K^{\prime}}$ for $z$$\in$$\Omega$. Then
$u(S(\varphi))$$H(\varphi)$=$P_{H(\varphi)}$$u(S_{K^{\prime}})H(\varphi)\subset{P_{H(\varphi)}}$$uH^{2}($$\Omega,$$K^{\prime}$$)
\subset{P_{H(\varphi)}}$$\varphi$
$H^{2}($$\Omega,$$K$$)$. Thus $u(S(\varphi))$=$0$. Thus
$(a)$$\rightarrow$$ (b)$.

Next, assume $(b)$. Then
$u(S_{K^{\prime}})H(\varphi)$=$uH(\varphi)$$\subset$$\varphi$$H^{2}($$\Omega,$$K$$)$.
It follows that $uH^{2}($$\Omega,$
$K^{\prime}$$)$
=$uH(\varphi)$$+$$u\varphi$$H^{2}($$\Omega,$$K$$)$$\subset$$\varphi$$H^{2}($$\Omega,$$K$$)$.
Thus $(b)$$ \rightarrow$$ (c)$.

To prove $(c)$$ \rightarrow$$(a)$, let
$M$=$\{$$f$$\in$$H^{2}($$\Omega,$$K$$)$$:$$ug$=$\varphi$$f$ for
some $g$$\in$$H^{2}($$\Omega,$$K^{\prime}$$)$$\}$. Then
$\overline{M}$ is $R(\Omega)$-invariant for $S_{K}$. By Theorem
1.4.2, there is a Hilbert space $K_{1}$ and a quasi-inner function
$\varphi_{1}$$\in$${H}^{\infty}$$(\Omega$,$L(K_{1}$,$K$$)$$)$ such
that $\overline{M}$=$\varphi_{1}$$H^{2}($$\Omega,$$K_{1}$$)$. From
Theorem 2.2.4 in [20], $u$=$\theta$$F$ where $\theta$ is a
function such that $\mid$$\theta$$\mid$ is constant almost
everywhere on each component of $\partial\Omega$ and $F$ is an
outer function in ${H}^{\infty}$$(\Omega$$)$. By the definition of
$M$,
${{\theta}H^{2}(\Omega,K^{\prime})}$=$\overline{{\theta}FH^{2}(\Omega,K^{\prime})}$=$\overline{uH^{2}(\Omega,K^{\prime})}$=$\overline{{\varphi}M}$=$\varphi$$\overline{M}$=${\varphi}{\varphi}_{1}$$H^{2}($$\Omega,$$K_{1}$$)$.
Since $\theta$ is quasi-inner,
${\theta}I_{K^{\prime}}$$\in$${H}^{\infty}$$(\Omega$,$L(K^{\prime}))$
is quasi-inner.(Note that
$({\theta}I_{K^{\prime}})(z)$=$\theta(z)I_{K^{\prime}}$). Then by
Corollary 1.4.4, there exist
$\varphi_{2}$$\in$${H}^{\infty}$$(\Omega$,$L(K^{\prime}$,$K_{1}$$)$$)$
such that
${\theta}I_{K^{\prime}}$=$\varphi\varphi_{1}\varphi_{2}$. Then
$uI_{K^{\prime}}$=$\varphi(F\varphi_{1}\varphi_{2})$, i.e.
$u(z)I_{K^{\prime}}$=$\varphi(z)(F(z)\varphi_{1}(z)\varphi_{2}(z))$.
Since $F\varphi_{1}\varphi_{2}$$\in$${H}^{\infty}$$(\Omega$,
$L(K^{\prime}$,$K$$)$$)$, $u$ is a scalar multiple of $\varphi$.
\end{proof}
In the next statement, $\adj
\varphi:\Omega\rightarrow{L({\C}^{n})}$ is defined by ($\adj
\varphi)(z)$=$\adj(\varphi(z))$ which is the algebraic adjoint of
$\varphi(z)$(i.e.,
$A\cdot\adj(A)$=$\adj(A)\cdot{A}$=$\det(A)I_{{\C}^{n}}$ for
$A\in{L({\C}^{n})})$.
\begin{prop}
Let $K$ and $K^{\prime}$ be Hilbert spaces with $\dim$ $K=$ $\dim$
$K^{\prime}=n(<\infty)$.

(a) If $\varphi$$\in$$H^{\infty}($$\Omega,$$L({\C}^{n})$$)$ is a
quasi-inner function, then $\theta$, defined by $\theta(z)$=

$\det$$(\varphi(z))$, is quasi-inner.

(b) If $\varphi$$\in$$H^{\infty}($$\Omega,$$L({\C}^{n})$$)$ is a
quasi-inner function, then $\adj\varphi$ is quasi-inner.

(c) If
$\varphi$$\in$${H}^{\infty}$$(\Omega$,$L(K$,$K^{\prime}$$)$$)$ is
a quasi-inner function, then $S(\varphi)$ is of class $C_{0}$.
\end{prop}

\begin {proof}
(a) and (b):  Since $\varphi$ is quasi-inner, there exists $m>0$
such that  for $h$$\in$${\C}^{n}$,
$m\norm{h}\leq\norm{\varphi(z)h}$ a.e. $z\in\partial\Omega$. Then
$m^{n}\leq$$\mid$$\det(\varphi(z))$$\mid$ and
$m^{n-1}\leq$$\mid$$\adj(\varphi(z))$$\mid$ a.e.
$z\in\partial\Omega$. From those facts one can conclude that (a)
and (b) are true.

(c): By Theorem 1.7.2, it is enough to prove that $\varphi$ has a
scalar multiple $u$$\in$${H}^{\infty}$$(\Omega)$. Let
$\psi(z)$=$\adj(\varphi(z))$ and $u(z)$=$\det(\varphi(z))$. Then
by (b), $\psi$$\in$${H}^{\infty}$$(\Omega$,$L(K^{\prime}$,
$K$$)$$)$ and by (a), $u(\neq{0})$$\in$${H}^{\infty}$$(\Omega)$.
Since $\varphi(z)[\adj
(\varphi(z))]$=$[\det(\varphi(z))]I_{K^{\prime}}$ for
$z\in\Omega$, it is proven.
\end{proof}
Let $\theta$ and $\theta^{\prime}$ be two  quasi-inner functions
in $H^{\infty}(\Omega)$. If $\theta\equiv{\theta}^{\prime}$ i.e.,
$\theta$ and $\theta^{\prime}$ belong to the same equivalence
class under the equivalence relation $\equiv$ between
$H^{\infty}(\Omega)$ functions introduced after Definition 1.6.1.,
then it is convenient to regard them as the same element in
$H^{\infty}(\Omega)$, and introduce the following definition.
\begin{defn}
Let $F$ be a family of functions in $H^{\infty}(\Omega)$. A
quasi-inner function $\theta$$\in$$H^{\infty}(\Omega)$ is called
\emph{the greatest common quasi-inner divisor} of $F$ if $\theta$
divides every element in $F$ and if $\theta$ is a multiple of any
other common quasi-inner divisor of $F$. \emph{The greatest common
quasi-inner divisor} of $F$ is denoted by $\bigwedge$$F$$($or
$\bigwedge$$_{i\in{I}}$$f_{i}$ if $F$=$\{f_{i}:$ $i$$\in$$I$ $\}$,
or $f_{1}\wedge$$f_{2}$ if $F$=$\{f_{1},$ $f_{2}\})$.
\end{defn}

\section{Quasi-equivalence and Quasi-similarity}
\goodbreak
\vskip0.1cm
\subsection{Normal form }
\begin{defn}
A \emph{quasi-unit} $\textbf{X}$ of order $n$ is a collection of
$n\times{n}$ matrices over $H^{\infty}(\Omega)$ such that the
family $\det(\textbf{X})$=$\{\det(X):$ $X$$\in$$\textbf{X}\}$ is
relatively prime, i.e. $\bigwedge\det(\textbf{X})\equiv{1}$.
\end{defn}
\begin{defn}
If $A$ and $B$ are $m\times{n}$ matrices over
$H^{\infty}(\Omega)$, then $A$ is said to be
\emph{quasi-equivalent} to $B$ if there exist quasi-units
$\textbf{X}$ and $\textbf{Y}$ of order $m$ and $n$ respectively
such that $\textbf{X}A$=${B}\textbf{Y}$ where
$\textbf{X}A$=$\{XA:X\in{\textbf{X}}\}$ and
$B\textbf{Y}$=$\{BY:Y\in{\textbf{Y}}\}$.
\end{defn}

A matrix $E$ over $H^{\infty}(\Omega)$ is in \emph{normal form}
$($or simply \emph{normal}$)$ provided
\begin{equation}E=\begin{pmatrix} D&0\\0&0
\end{pmatrix}\end{equation}
where $D$ is a diagonal matrix of nonzero quasi-inner functions
and each one except the first divides its predecessor.
\begin{defn}
Let $D_{k}(A)$ be the greatest common quasi-inner divisor of all
minors of rank $k$ of $A$ $($$k$ is no larger than min$\{m,n\})$
and $D_{0}$=$1$. Then the \emph{invariant factors} for a
$m\times{n}$ matrix $A$ over $H^{\infty}(\Omega)$ are defined by
$\xi_{k}(A)$=$D_{k}(A)/D_{k-1}(A)$ for $k$$\geq$$1$ such that some
minors of rank $k$ are not $0$.
\end{defn}

The following result is proved as Theorem 3.1 in [14].
\begin{prop}
Every $n$$\times{n}$ matrix over $H^{\infty}(\Omega)$ is
quasi-equivalent to a normal matrix. In fact, for any
$n$$\times{n}$ matrix $A$ over $H^{\infty}(\Omega)$, $A$ is
quasi-equivalent to the normal matrix formed by the invariant
factors of $A$.
\end{prop}

The following result is proved as in the case of the open unit
disk [6].
\begin{cor}
Let $\varphi$ be a quasi-inner function in
$H^{\infty}($$\Omega$,$L({\C}^{n}))$. If $\varphi$ is
quasi-equivalent to a normal matrix $N$ whose diagonal entries are
$\theta_{0},\dots,\theta_{n-1}$, then $\det$
$\varphi\equiv\theta_{0}\cdots\theta_{n-1}$.
\end{cor}
Let $f_{1}$ and $f_{2}$ be in $H^{\infty}(\Omega)$. If $M$ is the
$w^{*}$-closure of $f_{1}$$H^{\infty}(\Omega)$
$+$$f_{2}$$H^{\infty}(\Omega)$, then by the same way as Theorem 1
in [15], we can get $M$=$(f_{1}\wedge$$f_{2})H^{\infty}(\Omega)$.
\begin{prop}
Let $\varphi_{1}$, $\varphi_{2}$$\in$$H^{\infty}($$\Omega)$ be
such that $\varphi_{1}$$\wedge$$\varphi_{2}$ $\equiv$ $1$. If
$f$$\in$$L^{2}(\partial\Omega$,${\C}^{n})$ and $\varphi_{1}f$,
$\varphi_{2}f$$\in$$H^{2}(\partial\Omega$,${\C}^{n})$, then
$f$$\in$$H^{2}(\partial\Omega$,${\C}^{n})$.
\end{prop}
\begin{proof}
Since $\varphi_{1}$$\wedge$$\varphi_{2}\equiv{1}$, the
$w^{*}$-closure of
$\varphi_{1}$$H^{\infty}(\partial\Omega)$$+$$\varphi_{2}$$H^{\infty}(\partial\Omega)$
is $H^{\infty}(\partial\Omega)$. Thus there are nets
$\{f_{\alpha}\}$ and $\{g_{\alpha}\}$ in
$H^{\infty}(\partial\Omega)$ such that
$h_{\alpha}$=$\varphi_{1}$$f_{\alpha}$$+$$\varphi_{2}$$g_{\alpha}$
converges to 1, i.e.
\begin{equation}\int_{\partial\Omega} (h_{\alpha}-1)h dm
\rightarrow 0 \end{equation}for any
$h$$\in$$L^{1}(\partial\Omega)$. We will prove that
$h_{\alpha}f$$\rightarrow$$f$ weakly in $L^{2}(\partial\Omega$,
${\C}^{n})$, i.e. $((h_{\alpha}f$-$f)$, $g)$$\rightarrow$$0$ for
any $g$$\in$$L^{2}(\partial\Omega$,${\C}^{n})$. Indeed, if
$f$=$(f_{1},\dots,{f_{n}})$ and $g$=$(g_{1},\dots,{g_{n}})$, then
$(h_{\alpha}f-f$, $g)=\int_{\partial\Omega}$
$(h_{\alpha}$$-$$1)hdm$ where
$h=\sum_{i}f_{i}\bar{g_{i}}\in{L^{1}(\partial\Omega)}$. From
(2.1.2), we have \[h_{\alpha}f\rightarrow{f}\texttt{ weakly in
 } L^{2}(\partial\Omega).\] Since a subspace of a Banach space is
norm closed if and only if it is weakly closed [9],
$H^{2}(\partial\Omega$,${\C}^{n})$ is weakly closed. Since
$\varphi_{1}f$,
$\varphi_{2}f$$\in$$H^{2}(\partial\Omega$,${\C}^{n})$,
$h_{\alpha}f$$\in$$H^{2}(\partial\Omega$,${\C}^{n})$. If follows
that $f$$\in$$H^{2}(\partial\Omega$,${\C}^{n})$.
\end{proof}
The following results are proved as in the case of the open unit
disk [6].
\begin{prop}
Let $\varphi_{1}$ and $\varphi_{2}$ be quasi-inner functions in
$H^{\infty}($$\Omega$,$L({\C}^{n}))$. If $\varphi_{1}$ and
$\varphi_{2}$ are quasi-equivalent, then $S(\varphi_{1})$ and
$S(\varphi_{2})$ are quasisimilar.
\end{prop}
\begin{cor} Let $\varphi$ be a quasi-inner function in $H^{\infty}($$\Omega$,$L({\C}^{n}))$. If $\varphi$ is
quasi-equivalent to a normal matrix $N$ whose diagonal entries are
$\theta_{0},\dots,\theta_{n-1} ($$\theta_{i+1}$$\mid$$\theta_{i}$
for $i$=$0,1,\dots,{n-1})$, then $S(\varphi)\sim{
\oplus_{i=0}^{n-1}S(\theta_{i})}$.\end{cor}
\begin{proof} Since $S(N)\sim{\bigoplus_{i=0}^{n-1}S(\theta_{i})}$, by Proposition
2.1.7, $S(\varphi)\sim{\bigoplus_{i=0}^{n-1}S(\theta_{i})}$,
because $"\sim"$ is an equivalence relation.\end{proof}
\begin{cor}Let $\varphi_{1}$ and $\varphi_{2}$ be quasi-inner
functions in $H^{\infty}($$\Omega$,$L({\C}^{n}))$. If
$S(\varphi_{1})$ is a quasi-affine transform of $S(\varphi_{2})$,
then $\varphi_{1}$ and $\varphi_{2}$ are
quasi-equivalent.\end{cor}
\begin{cor}Let $\varphi$  be a
quasi-inner function in $H^{\infty}($$\Omega$,$L({\C}^{n}))$. If
$S(\varphi)\sim{\bigoplus_{i=0}^{n-1}S(\theta_{i})}$, then $det$
$\varphi\equiv\theta_{0}\cdots\theta_{n-1}$.\end{cor}\begin{proof}
Let $N$ be a normal matrix whose diagonal entries are
$\theta_{0},\dots,\theta_{n-1}$. Since $S(N)\sim{
\bigoplus_{i=0}^{n-1}S(\theta_{i})}$, $S(\varphi)\sim{S(N)}$. By
Corollary 2.1.9, $\varphi$ and $N$ are quasi-equivalent. Then by
Corollary 2.1.5, det
$\varphi\equiv\theta_{0}\cdots\theta_{n-1}$\end{proof}
\subsection{Main results}
In this section, first of all we show how to use Theorem 1.4.2 and
Corollary 1.4.4.

\begin{thm}Let $F$ and $F^{\prime}$ be two separable Hilbert spaces and $\varphi$  be a
quasi-inner function in $H^{\infty}($$\Omega$,$L(F,F^{\prime}))$.

$(i)$ If $M$ $\subset$ $H(\varphi)$ is $R(\Omega)$-invariant for
$S(\varphi)$, then there is a Hilbert space $K$ and there are
quasi-inner functions $\varphi_{1}$ $\in$
$H^{\infty}($$\Omega$,$L(F,K))$ and
$\varphi_{2}$$\in$$H^{\infty}($$\Omega$,$L(K,F^{\prime}))$ such
that $\varphi(z)$=$\varphi_{2}(z)\varphi_{1}(z)$ for
$z$$\in$$\Omega$ and
\begin{equation}M=\varphi_{2}H^{2}(\Omega,K)\ominus\varphi{H^{2}(\Omega,F)}\end{equation}

$(ii)$ Conversely, if $K$, $\varphi_{1}$ and $\varphi_{2}$ are as
above , then $(2.2.1)$ defines a $R(\Omega)$-invariant subspace of
$H(\varphi)$. Moreover, if $S(\varphi)$=$\begin{pmatrix}
T_{1}&X\\0&T_{2}\end{pmatrix}$ is the triangularization of
$S(\varphi)$ with respect to the decomposition
$H(\varphi)$=$M\oplus(H(\varphi)$ $\ominus{M})$, then
$T_{2}$=$S(\varphi_{2})$ and $S(\varphi_{1})$ is similar to
$T_{1}$.
\end{thm}
\begin{proof} $(i).$
Since $M$ is  $R(\Omega)$-invariant, the space
$M\oplus\varphi{H^{2}(\Omega,F)}$ is also $R(\Omega)$-invariant
subspace of ${H^{2}(\Omega,F^{\prime})}$ and so Theorem 1.4.2
implies the existence of a Hilbert space $K$ and of a quasi-inner
function
$\varphi_{2}$$\in$$H^{\infty}($$\Omega$,$L(K,F^{\prime}))$ such
that $(2.2.1)$ holds.

The inclusion
$\varphi{H^{2}(\Omega,F)}\subset\varphi_{2}{H^{2}(\Omega,K)}$
implies that for any $f$$\in$${H^{2}(\Omega,F)}$ there is
$\phi_{f}$$\in$${H^{2}(\Omega,K)}$ such that
$\varphi{f}$=$\varphi_{2}\phi_{f}$. Let
$M^{\prime}$=$\{$$\phi_{f}$$\in$${H^{2}(\Omega,K)}:$ $\varphi{f}$
=$\varphi_{2}\phi_{f}$ for some $f$$\in$${H^{2}(\Omega,F)}$$\}$.
Since $\varphi(rf)$=$\varphi_{2}(r\phi_{f})$ for any
$r$$\in$$R(\Omega)$ and $f$$\in$${H^{2}(\Omega,F)}$, $M^{\prime}$
is also a $R(\Omega)$-invariant subspace of ${H^{2}(\Omega,K)}$,
and so $M^{\prime}$=$\varphi_{3}{H^{2}(\Omega,K^{\prime})}$ for
some Hilbert space $K^{\prime}$ and a quasi-inner function
$\varphi_{3}$$\in$$H^{\infty}($$\Omega$,$L(K^{\prime},K))$ by
theorem 1.4.2 It follows that
$\varphi{H^{2}(\Omega,F)}$=$\varphi_{2}$$\varphi_{3}{H^{2}(\Omega,K^{\prime})}$
by the definition of $M^{\prime}$. By Corollary 1.4.4, there is a
function
$\varphi_{4}$$\in$$H^{\infty}($$\Omega$,$L(F,K^{\prime}))$ such
that $\varphi$=$\varphi_{2}\varphi_{3}\varphi_{4}$.

Let
$\varphi_{1}$=$\varphi_{3}\varphi_{4}$$\in$$H^{\infty}($$\Omega$,$L(F,K))$.
Since $\varphi$ and $\varphi_{2}$ are quasi-inner functions, so is
$\varphi_{1}$. Thus $\varphi_{1}$ is a quasi-inner function
satisfying $\varphi$=$\varphi_{2}\varphi_{1}$.\vskip0.3cm

$(ii).$ The $R(\Omega)$-invariance of the subspace $M$ described
by $(2.2.1)$ is obvious.

Since
$H(\varphi)\ominus{M}$=${H^{2}(\Omega,F^{\prime})}\ominus\varphi_{2}{H^{2}(\Omega,K)}$=$H(\varphi_{2})$,
we have
${T_{2}}^{*}$=${S(\varphi)}^{*}$$|$$H(\varphi)\ominus{M}$=$S^{*}_{F^{\prime}}$$|$$H(\varphi_{2})$=${S(\varphi_{2})}^{*}$.
Thus $T_{2}$=$S(\varphi_{2})$. It remains to prove similarity of
$T_{1}$ and $S(\varphi_{1})$. Define
$Y:H^{2}(\Omega,K)\rightarrow\varphi_{2}$$H^{2}(\Omega,K)$ by
$Yf=\varphi_{2}f$. Clearly $Y$ is onto. Since $\varphi_{2}$ is a
quasi-inner function, $Y$ is one-to-one. Since
$Y(\varphi_{1}$$H^{2}(\Omega,F))$=$\varphi_{2}\varphi_{1}$$H^{2}(\Omega,F)$
=$\varphi$$H^{2}(\Omega,F)$,
$\varphi_{2}$$H^{2}(\Omega,K)$=$M\oplus\varphi$$H^{2}(\Omega,F)$
and
$H^{2}(\Omega,K)$=$H(\varphi_{1})\oplus\varphi_{1}$$H^{2}(\Omega,F)$,
we have $P_{M}Y(H(\varphi_{1}))$=$M$. Thus we can define a bounded
linear operator $F:H(\varphi_{1})\rightarrow{M}$ by
$Fg=P_{M}$$\varphi_{2}g$ for $g\in{H(\varphi_{1})}$, and $F$ is
onto. Since $\varphi_{2}$ is a quasi-inner function, ker
$F=\{g\in{H(\varphi_{1})}$$:$
$\varphi_{2}g\in{\varphi{H^{2}(\Omega,F)\}}}$
=$\{g\in{H(\varphi_{1})}:g$$\in\varphi_{1}$$H^{2}(\Omega,F)\}=\{0\}$.
It follows that $F\in{L(H(\varphi_{1}),M)}$ is bijective. By the
Open Mapping Theorem, $F$ is invertible and clearly
$T_{1}F=FS(\varphi_{1})$.
\end{proof}
\vskip0.3cm

Fix $n\geq{1}$, and consider the mapping
$\Gamma_{n}:L(F)\rightarrow{L({\otimes}^{n}F)}$ given by
$\Gamma_{n}(T)=T\otimes{T}\otimes\cdot\cdot\cdot\otimes{T}$, where
$F$ is a Hilbert space and $T\in{L(F)}$.

Define a unitary representation
$\pi_{n}:S_{n}\rightarrow{L({\otimes}^{n}F)}$ where $S_{n}$
denotes the group of permutations of $\{1,2,\cdot\cdot\cdot,n\}$,
defined by

\begin{equation}\pi_{n}(\sigma)(x_{1}\otimes{x_{2}}\otimes\cdot\cdot\cdot\otimes{x_{n}})=x_{\sigma^{-1}(1)}\otimes{x_{\sigma^{-1}(2)}}\otimes\cdot\cdot\cdot\otimes{x_{\sigma^{-1}(n)}},\end{equation}
$\sigma\in{S_{n}}$, $x_{j}\in{F}$, $1\leq{j}\leq{n}$.

The homomorphism $\pi_{n}$ can be extended to a $*$-homomorphism,
still denoted $\pi_{n}$, from the $C^{*}$-algebra consisting of
all formal sums $\sum_{\sigma\in{S_{n}}}\alpha_{\sigma}\sigma$
$(\alpha_{\sigma}\in{\C})$ to $L({\otimes}^{n}F)$. We will use the
alternating projection $a_{n}$ defined by
\begin{equation}a_{n}=\frac{1}{n!}\sum_{\sigma\in{S_{n}}}\epsilon(\sigma)\sigma,\end{equation}
where $\epsilon(\sigma)$ is the sign of $\sigma$, i.e.
$\epsilon(\sigma)=+1$ or $-1$ according to whether $\sigma$ is an
even or odd permutation. Let $n\geq{1}$ be a natural number.We use
the notation $\bigwedge^{n}$${F}$ for
$\pi_{n}(a_{n})({\otimes}^{n}F)$. The space $\bigwedge^{n}$${F}$
is called the \emph{n}th \emph{exterior power} of $F$. If
$B\in{L(F)}$, we denote by $\bigwedge^{n}$${B}$ the operator
$\Gamma_{n}(B)|\bigwedge^{n}$${F}$.
\begin{prop}
If $A$ and $B$ are quasi-equivalent quasi-inner functions in
$H^{\infty}($$\Omega$, $L({\C}^{n}))$, then $\bigwedge^{k}$${A}$
and $\bigwedge^{k}$${B}$ are quasi-equivalent, for
$1\leq{k}\leq{n}$.
\end{prop}
\begin{proof}
This is same as Proposition 6.5.17 in [6].
\end{proof}
\begin{prop}
If $A=\begin{pmatrix}
\theta_{0}&0&\dots&{0}\\0&\theta_{1}&\dots&{0}\\\\\\{0}&0&\dots&\theta_{n-1}\end{pmatrix}_{n\times{n}}$
is normal, then $S(\bigwedge^{k}$${A})$ has minimal function
$\theta_{0}\theta_{1}\cdots\theta_{k-1}$ for $k=1,\dots,{n}$.
\end{prop}
\begin{proof}
 Since
$\bigwedge^{k}$${A}$ is also a diagonal quasi-inner function with
diagonal entries
$\theta_{i_{1}}\theta_{i_{2}}\cdot\cdot\cdot\theta_{i_{k}}$ where
$i_{p}\neq{i_{q}}$ for $p\neq{q}$ ([6]), the minimal function of
$S(\bigwedge^{k}$${A})$ is
$\theta_{0}\theta_{1}\cdot\cdot\cdot\theta_{k-1}$.
\end{proof}
If $\{M_{i}\}_{i\in{I}}$ is a family of subsets of the Hilbert
space $H$, we denote by $\bigvee_{i\in{I}}M_{i}$ the closed linear
span generated by $\bigcup_{i\in{I}}M_{i}$.

\begin{defn}
Let $T$ $\in$ ${L(H)}$ be an operator with spectrum in
$\overline{\Omega}$. A subset $G$ $\subseteq$ $H$ with the
property that $\bigvee$$\{r(T)m$ $;$ $r$ $\in$ $R(\Omega),$ $m$
$\in$ $G\}$ = $H$, is called an $R(\Omega)$-\emph{cyclic set} for
$T$. The \emph{multiplicity} $\mu_{T}$ of $T$ is the smallest
cardinality of an $R(\Omega)$-cyclic set for $T$. The operator $T$
is said to be \emph{multiplicity-free} if $\mu_{T}$ = $1$. If
$\mu_{T}$ = $1$, any vector $x$ $\in$ $H$ such that
$\bigvee$$\{r(T)x$ $;$ $r$ $\in$ $R(\Omega)\}$ = $H$ is said to be
$R(\Omega)$-\emph{cyclic} for $T$.
\end{defn}

Recall that if $\mu_{T}\leq{n}$, then Jordan model of $T$ is
$\oplus_{j=0}^{n-1}S(\theta_{j})$ \cite{20}.
\begin{prop}
Assume that $T$ $\in$ $L(H)$ is an operator of class $C_{0}$
relative to ${\Omega}$ such that $\mu_{T}$ = $n$ $<$ $\infty$,
$H^{\prime}$ is a $R(\Omega)$-invariant subspace for $T$, and $T$
= $\begin{pmatrix} T^{\prime}&Y\\0&T^{\prime\prime}\end{pmatrix}$
is the triangularization of $T$ with respect to the decomposition
$H$ = $H^{\prime}$ $\oplus$ $(H\ominus{H^{\prime}})$. If
$\oplus_{j<n}S(\theta_{j})$, $\oplus_{j<n}S(\theta_{j}^{\prime})$,
and $\oplus_{j<n}S(\theta_{j}^{\prime\prime})$ are the Jordan
models of $T$, $T^{\prime}$, $T^{\prime\prime}$, respectively,
then
\[\theta_{0}\cdot\cdot\cdot\theta_{k-1}|\theta_{0}^{\prime}\cdot\cdot\cdot\theta_{k-1}^{\prime}
\theta_{0}^{\prime\prime}\cdot\cdot\cdot\theta_{k-1}^{\prime\prime}\]
for every $k$ such that $1$$\leq$$k$$<$$n$, and
\[\theta_{0}\cdot\cdot\cdot\theta_{n-1}\equiv\theta_{0}^{\prime}
\cdot\cdot\cdot\theta_{n-1}^{\prime}\theta_{0}^{\prime\prime}\cdot\cdot\cdot\theta_{n-1}^{\prime\prime}.\]
\end{prop}
\begin{proof}
Let $f$ $\in$ $H^{\infty}($$\Omega$,$L({\C}^{n}))$ be a
quasi-inner function such that $f$ is a normal matrix whose
diagonal entries are $\theta_{0},\dots,\theta_{n-1}$. By Corollary
2.1.8, $S(f)=\oplus_{j=0}^{n-1}S(\theta_{j})\sim{T}$. Thus there
is an injective operator $X\in{L(H,H(f))}$ with dense range such
that $S(f)X=XT$.

Let $M$ be the closure of $XH^{\prime}$. Since $H^{\prime}$ is a
$R(\Omega)$-invariant subspace for $T$, so is $M$ for $S(f)$. Then
by Theorem 2.2.1, there are quasi-inner functions
$f_{1}\in{H^{\infty}(\Omega,L({\C}^{n}))}$ and
$f_{2}\in{H^{\infty}(\Omega,L({\C}^{n}))}$ such that
$f=f_{2}f_{1}$ and
$M=f_{2}H^{2}(\Omega,$${\C}^{n})\ominus{f}H^{2}(\Omega,$${\C}^{n})$.
If $S(f)=\begin{pmatrix} T_{1}&Z\\0&T_{2}\end{pmatrix}$ is the
triangularization of $S(f)$ with respect to the decomposition
$H(f)=M\oplus{(H(f)\ominus{M})}$, then by Theorem 2.2.1, $T_{1}$
is similar to $S(f_{1})$ and $T_{2}=S(f_{2})$.

Let $X^{\prime}=X$$|$$H^{\prime}$. Then
$T_{1}X^{\prime}=S(f)X$$|$$H^{\prime}=XT$$|$$H^{\prime}=X^{\prime}T^{\prime}$
and so $T_{1}\sim{T^{\prime}}\sim
\oplus_{j=0}^{n-1}S(\theta_{j}^{\prime})$. Since $T_{1}$ is
similar to $S(f_{1})$,
$S(f_{1})\sim\oplus_{j=0}^{n-1}S(\theta_{j}^{\prime})$. Define
$X^{\prime\prime}:H(f)\ominus{M}\rightarrow{H}\ominus{H^{\prime}}$
by $X^{\prime\prime}=X^{\ast}$$|$$H(f)\ominus{M}$. Then
$X^{\prime\prime}$ is injective with dense range in
$H\ominus{H^{\prime}}$ and $X^{\prime\prime}$${T_{2}}^{\ast}=
X^{\ast}S(f)^{\ast}$$|$$H(f)\ominus{M}=T^{\ast}X^{\ast}$$|$$H(f)\ominus{M=(T^{\prime\prime})^{\ast}X^{\prime\prime}}$.
Thus
$T_{2}\sim{T^{\prime\prime}}\sim\oplus_{j=0}^{n-1}S(\theta_{j}^{\prime\prime})$.
It follows that
$S(f_{2})\sim\oplus_{j=0}^{n-1}S(\theta_{j}^{\prime\prime})$. Fix
$k$ such that $1\leq{k}<{n}$ and note that
$\bigwedge^{k}$$f=\bigwedge^{k}f_{2}\wedge\bigwedge^{k}$$f_{1}$.
By Proposition 2.2.3, the minimal function of $S(\bigwedge^{k}f)$
is $\theta_{0}\theta_{1}\cdots\theta_{k-1}$. By Corollary 2.1.9,
there are normal matrices $N_{1}$ and $N_{2}$ which are
quasi-equivalent to $f_{1}$ and $f_{2}$, respectively and diagonal
entries of $N_{1}$ ($N_{2}$) are
$\theta^{\prime}_{0},\theta^{\prime}_{1},\dots,\theta^{\prime}_{n-1}$
($\theta^{\prime\prime}_{0},\theta^{\prime\prime}_{1},\dots,\theta^{\prime\prime}_{n-1}$
, respectively). By Proposition 2.2.2, $\bigwedge^{k}f_{1}$ and
$\bigwedge^{k}N_{1}$ are quasi-equivalent. By Proposition 2.1.7,
$S(\bigwedge^{k}f_{1})$ and $S(\bigwedge^{k}N_{1})$ are
quasisimilar. Thus the minimal functions of
$S(\bigwedge^{k}f_{1})$ is
$\theta_{0}^{\prime}\theta_{1}^{\prime}\cdots\theta_{k-1}^{\prime}$.
Similarly, the minimal function of $S(\bigwedge^{k}f_{2})$ is
$\theta_{0}^{\prime\prime}\theta_{1}^{\prime\prime}\cdots\theta_{k-1}^{\prime\prime}$.
By Theorem 1.7.2, there are functions $g^{\prime}$,
$g^{\prime\prime}\in{H^{\infty}(\Omega,L(\bigwedge^{k}{\C}^{n}))}$
such that
$g^{\prime}(\bigwedge^{k}f_{1})=\theta_{0}^{\prime}\theta_{1}^{\prime}\cdots\theta_{k-1}^{\prime}$$I$
and
$g^{\prime\prime}(\bigwedge^{k}f_{2})=\theta_{0}^{\prime\prime}\theta_{1}^{\prime\prime}\cdots\theta_{k-1}^{\prime\prime}$$I$.
Combining these relations we get
$g^{\prime}g^{\prime\prime}(\bigwedge^{k}f)=g^{\prime}g^{\prime\prime}(\bigwedge^{k}f_{2}\bigwedge^{k}f_{1})=
\theta_{0}^{\prime}\theta_{1}^{\prime}\cdots\theta_{k-1}^{\prime}\theta_{0}^{\prime\prime}\theta_{1}^{\prime\prime}\cdots\theta_{k-1}^{\prime\prime}I$
and this Corollary follows because
$\theta_{0}\theta_{1}\cdots\theta_{k-1}$ is the least scalar
multiple of $\bigwedge^{k}f$ by Theorem 1.7.2.

Next, for $k=n$, since $S(f)\sim\oplus_{j=0}^{n-1}S(\theta_{j})$,
$S(f_{1})\sim\oplus_{j=0}^{n-1}S(\theta_{j}^{\prime})$, and
$S(f_{2})\sim\oplus_{j=0}^{n-1}S(\theta_{j}^{\prime\prime})$, by
Corollary 2.1.10,
$\det(f)\equiv\theta_{0}\theta_{1}\cdots\theta_{n-1}$,
$\det(f_{1})\equiv\theta_{0}^{\prime}\theta_{1}^{\prime}\cdots\theta_{n-1}^{\prime}$,
and
det$(f_{2})\equiv\theta_{0}^{\prime\prime}\theta_{1}^{\prime\prime}\cdots\theta_{n-1}^{\prime\prime}$.
From the fact $f=f_{2}f_{1}$, we can get
det$(f)=($det$(f_{2}))($det$(f_{1}))$ which proves the case $k=n$.
\end{proof}
When $T\in{L(H)}$ is an operator of class $C_{0}$ relative to
${\Omega}$ and $K=\bigvee\{r(T)h:$ $r\in{R(\Omega)}\}$, let $m_h$
denote the minimal function of $T|{K}$. We have the following
Proposition from Theorem 4.3.10. in [20].
\begin{prop}\label{k2}
Let $T\in{L(H)}$ be an operator of class $C_{0}$ relative to
${\Omega}$. If ${\oplus_{j<\omega}S(\theta_{j})}$ is the Jordan
model of $T$, then for any $k=1,2,3,...$, there are
$R(\Omega)$-invariant subspaces $M_{-1}$, $M_{0}$,..., $M_{k-2}$
and $h_{0}$, $h_{1}$,...,$h_{k-1}$ in $H$ such that

\begin{equation}\label{k3} h_{i}\in{M_{i-1}}\emph{ }and\texttt{ }
m_{h_{i}}=m_{T\mid{M_{i-1}}}\end{equation} for $i=0,1,...,k-1,$
and
\begin{equation}\label{k4} K_{i}\vee{M}_{i}=M_{i-1}\emph{ }and\emph{ }
K_{i}\cap{M_{i}}=\{0\}\end{equation} for $i=0,1,...,k-1,$ where
$M_{-1}=H$ and $K_{i}=\bigvee\{r(T)h_{i}:$
$r\in{R(\Omega)}\}$\end{prop}

\begin{thm}
Assume that $T\in{L(H)}$ is an operator of class $C_{0}$ relative
to ${\Omega}$, $H^{\prime}$ is a $R(\Omega)$-invariant subspace
for $T$, and $T=\begin{pmatrix}
T^{\prime}&Y\\0&T^{\prime\prime}\end{pmatrix}$ is the
triangularization of $T$ with respect to the decomposition
$H=H^{\prime}\oplus{(H\ominus{H^{\prime}})}$. If
$\oplus_{j<\gamma}S(\theta_{j})$,
$\oplus_{j<\gamma}S(\theta_{j}^{\prime})$, and
$\oplus_{j<\gamma}S(\theta_{j}^{\prime\prime})$ are the Jordan
models of $T$, $T^{\prime}$, $T^{\prime\prime}$, respectively,
then
$\theta_{0}\cdots\theta_{k-1}$$|$$\theta_{0}^{\prime}\cdots\theta_{k-1}^{\prime}\theta_{0}^{\prime\prime}\cdots\theta_{k-1}^{\prime\prime}$
for every $k=1,2,3,\dots$.
\end{thm}
\begin{proof}
Fix $k\geq{1}$. Since $T\sim{\oplus_{j<\omega}S(\theta_{j})}$, by
Proposition \ref{k2} and proof of Theorem 4.3.10 in [20], there is
a $R(\Omega)$-invariant subspace $M$ for $T$ such that
$T_{1}$$=T|M\sim\oplus_{j<k}S(\theta_{j})$. Clearly
${H_{1}}^{\prime}=M\cap{H}^{\prime}$ is $R(\Omega)-$invariant for
$T_{1}$.

Let $T_{1}=\begin{pmatrix}
T_{1}^{\prime}&{Y_{1}}\\0&T_{1}^{\prime\prime}\end{pmatrix}$ be
the triangularization of $T_{1}$ with respect to the decomposition
$M=H_{1}^{\prime}\oplus{(M\ominus{H_{1}^{\prime}})}$. If
$\oplus_{j}S(\phi_{j}^{\prime})$ and
$\oplus_{j}S(\phi_{j}^{\prime\prime})$ are Jordan models of
$T_{1}^{\prime}$ and $T_{1}^{\prime\prime}$, respectively, then
\begin{equation} \theta_{0}\cdots\theta_{k-1}\equiv\phi_{0}^{\prime}\cdots\phi_{k-1}^{\prime}\phi_{0}^{\prime\prime}\cdots\phi_{k-1}^{\prime\prime}\end{equation}
by Proposition 2.2.5. $($Note $\mu_{T_{1}}\leq{k}$$)$. Since
${T}^{\prime}$$\mid$${H_{1}}^{\prime}={T_{1}}^{\prime}$, by
Proposition 1.6.3,
${\phi_{i}}^{\prime}$$\mid$${\theta_{i}}^{\prime}$ for
$i=0,\dots,k-1.$

Next, let $H_{1}^{\prime\prime}=M\ominus{H_{1}^{\prime}}$,
${H}^{\prime\prime}=H\ominus{H^{\prime}}$, and
$X:{H_{1}}^{\prime\prime}\rightarrow{H}^{\prime\prime}$ be
orthogonal projection. If $a\in{ker}$$X$, then
$a\in{H^{\prime}\cap(M\ominus{H_{1}^{\prime}})}\subset{{H^{\prime}\cap{M}}}={H_{1}}^{\prime}$.
Since $a\in{H_{1}^{\prime\prime}}(=M\ominus{H_{1}^{\prime}})$,
$a=0$. Thus $X$ is one-to-one. Moreover, $H^{\prime}$ is invariant
for $T$, and $H^{\prime\prime}$ is invariant for $T^{*}$. Thus
$T^{*}P_{H^{\prime\prime}}=P_{H^{\prime\prime}}T^{*}P_{H^{\prime\prime}}$
and so
$P_{H^{\prime\prime}}T=(T^{*}P_{H^{\prime\prime}})^{*}=(P_{H^{\prime\prime}}T^{*}P_{H^{\prime\prime}})^{*}=P_{H^{\prime\prime}}TP_{H^{\prime\prime}}=T^{\prime\prime}P_{H^{\prime\prime}}$.
Since
$P_{H^{\prime\prime}}T_{1}^{\prime\prime}=
P_{H^{\prime\prime}}P_{M\ominus{H_{1}^{\prime}}}T
\mid{M\ominus{H_{1}^{\prime}}}=P_{H^{\prime\prime}}T\mid{M\ominus{H_{1}^{\prime}}}$,
 $T^{\prime\prime}X=XT_{1}^{\prime\prime}$. Since $X$ is
one-to-one, $T_{1}^{\prime\prime}$ is quasi-similar to a
restriction of $T^{\prime\prime}$ to an invariant subspace and so
we can get
${\phi_{i}}^{\prime\prime}$$\mid$${\theta_{i}}^{\prime\prime}$ for
$i=0,...,k-1$. Thus from (2.2.6), we can conclude that
$\theta_{0}\cdot\cdot\cdot\theta_{k-1}$$|$$\theta_{0}^{\prime}\cdot\cdot\cdot\theta_{k-1}^{\prime}\theta_{0}^{\prime\prime}\cdot\cdot\cdot\theta_{k-1}^{\prime\prime}$.
\end{proof}

------------------------------------------------------------------------

\bibliographystyle{amsplain}
\bibliography{xbib}
\end{document}